\documentclass[11pt]{amsart}
   \newlength{\extramargin}
\usepackage{amsfonts, amsmath,latexsym,verbatim,amscd,amssymb, stmaryrd, wasysym, mathrsfs}

\newtheorem{theorem}{Theorem}
\newtheorem{corollary}[theorem]{Corollary}

\newtheorem{lemma}[theorem]{Lemma}
\newtheorem{proposition}[theorem]{Proposition}

\numberwithin{equation}{section}
\numberwithin{theorem}{section}

 \newcommand{\Fun}{
 \makebox{
 $\text{\scriptsize{I}} \hspace{-.3 mm} \text{\scriptsize{I}}$
 \hspace{-.215 in} \raisebox{.05 in}{$\smile$}
 \hspace{-.25 in} \raisebox{-.055 in}{$\frown$}}}

 \newcommand{\FunForm}{
 \makebox{$\mathcal{I} \hspace{-1.25 mm} \mathcal{I}$}}
 
\DeclareMathOperator*{\FP}{\text{FP}}
\DeclareMathOperator*{\Res}{\text{Res}}

\DeclareMathOperator{\dvol}{\text{dvol}}

\DeclareMathOperator{\Ind}{\text{Ind}}

\renewcommand{\mod}{\text{ mod }}

\newcommand{\dzero}[1]{\frac{\partial}{\partial #1}\biggr\rvert_{#1 =0}}
\newcommand{\Pff}{\sideset{^{R}}{}\int \mathrm{Pff}}
\newcommand{\tr}{\text{tr}}
\newcommand{\Rind}{\sideset{^{R}}{}\Ind}
\newcommand{\Rint}{\sideset{^{R}}{}\int}
\newcommand{\Hint}{\sideset{^{H}}{}\int}
\newcommand{\Rintd}{\sideset{^{R}}{_0^{\delta}}\int}
\newcommand{\Hintd}{\sideset{^{H}}{_0^{\delta}}\int}

\newcommand{\script}[1]{\textsc{#1}}
\newcommand{\df}[1]{\mathfrak{#1}}
\newcommand{\curly}[1]{\mathcal{#1}}

\newcommand{\lrpar}[1]{\left( #1 \right)}
\newcommand{\lrspar}[1]{\left[ #1 \right]}

\title[Renormalizing Curvature Integrals]{Renormalizing Curvature Integrals on Poincar\'e-Einstein manifolds}
\author{Pierre Albin}
\date{April 5, 2005}
\address{Department of Mathematics, Stanford University}
\email{pierre@math.stanford.edu}

\begin{document}

\begin{abstract}
After analyzing renormalization schemes on a Poincar\'e-Einstein manifold, we study the renormalized integrals of scalar Riemannian invariants.
The behavior of the renormalized volume is well-known, and we show any scalar Riemannian invariant renormalizes similarly.
We consider characteristic forms and their behavior under a variation of the Poincar\'e-Einstein structure, and obtain, from the renormalized integral of the Pfaffian, an extension of the Gauss-Bonnet theorem.
\end{abstract}

\maketitle

\section{Introduction}

Let $M$ be the interior of a compact manifold with boundary 
$\bar{M}$.  A boundary defining function 
(hereafter, a bdf), $x$, is a 
smooth nonnegative function on $\bar{M}$ that vanishes precisely at the boundary, 
with non-degenerate differential there.
A metric $g$ on $M$ is said to be conformally compact if 
there is a bdf $x$ such that $\bar{g}=x^{2}g$ 
extends to a nondegenerate metric on $\bar{M}$. Conformally compact Einstein metrics, also known as Poincar\'e-Einstein or PE metrics, have been the object of many recent studies by geometers and physicists. 

A conformally compact metric on the interior of a manifold determines a conformal class of metrics on the boundary but does not single out any particular metric within that class.
In analogy to the well-known relation between the hyperbolic geometry of the ball and the conformal geometry of the sphere, Fefferman and Graham \cite{Fefferman-Graham} proposed studying the conformal geometry of the boundary of a manifold through the geometry of a PE metric on the interior.

More recently, the string-theory community has been interested in these manifolds because of the role they play in the Anti-de-Sitter/Conformal field theory (AdS/CFT) correspondence conjectured by Maldacena.
Witten \cite{Witten} clarified the conjecture as an equality of partition functions (normalizing factors for the probabilities of states). Computing the partition function of the conformal field theory (in the massless case) involves evaluating the Einstein action. For an Einstein manifold, this is a multiple of the volume which, however, is not finite. Witten showed that the volume could be renormalized while preserving covariance. Hennington and Skenderis carried out explicit computations in low dimensions and verified that this notion of renormalized volume is consistent with computations from, for instance, superconformal SU(4) gauge theory. Later, Graham and Witten \cite{Graham-Witten}, carried out a similar analysis of the area of (minimal) submanifolds of PE manifolds.

A mathematical survey of these ideas appeared in \cite{Graham}. The definition of renormalized volume (see section $\mathcal{x}$\ref{sec:PE} below) involves the choice of a bdf. One can restrict the choice of bdf to a natural class -- the ``special" bdfs -- for which the metric near the boundary has the particularly nice form
\begin{equation}\label{IntroFG1}
       g = \frac{dx^2 + G_x}{x^2}. 
\end{equation}
The corresponding renormalized volume of an even dimensional PE manifold is independent of the choice of bdf from this smaller set. In odd-dimensions, even this restriction does not produce an invariant, and different choices of special bdfs yield different renormalized volumes producing the so-called ``conformal anomaly".

In sections $\mathcal{x}$\ref{sec:AH}-$\mathcal{x}$\ref{ScalInv}, we analyze the scalar Riemannian invariants of a PE manifold. These are all complete contractions of the curvature and its covariant derivatives. They include the volume and also the heat invariants. We show that one can define renormalized integrals of these invariants, and their
dependence on the choice of bdf is just like that of the volume in the following sense (cf. Theorem \ref{thm:PEinvs} below).

\begin{theorem}\label{Intro1}
Let (M,g) be a Poincar\'e-Einstein manifold, and $P$ a scalar Riemannian invariant of $g$. Using any given special bdf $x$, we may define the renormalized integral of $P$ over $M$, denoted by
\[ \Rint P \dvol_g .\]
If $M$ is even-dimensional, this renormalized integral is independent of the choice of special bdf. If $M$ is odd-dimensional, and $P$ is not integrable, its renormalized integral necessarily depends on the choice of special bdf. In this case, $P$ determines a ``residue" integral on $\partial \bar{M}$ which is independent of the choice of special bdf.
\end{theorem}

Of particular interest is the renormalized integral of the Pfaffian. Recall that the Pfaffian is a natural density that integrates to the Euler characteristic on closed even-dimensional manifolds. On an even-dimensional PE manifold the integral of the Pfaffian diverges, but the renormalization process recovers the same topological information as in the compact setting:

\begin{theorem}\label{Intro2}
On an even dimensional Poincar\'e-Einstein manifold, using a special bdf to renormalize the integral of the Pfaffian,
\begin{equation}\label{IntroPff}
        \Pff = \chi(M) .
\end{equation}
\end{theorem}

This generalizes the four-dimensional formula of Anderson \cite{Anderson},
\begin{equation}\label{IntroAnderson}
    \frac{1}{8(2\pi)^2}\int |W|^2 + \frac{3}{(2\pi)^2}\hat{V} = \chi(M) ,
\end{equation}
where $W$ is the Weyl curvature tensor and $\hat{V}$ is the renormalized volume. 
Notice that, by virtue of the Einstein condition, the Pfaffian is a polynomial in the Weyl curvature (see Lemma \ref{EinPfaff} for an explicit formula).
Another particular case is the formula of Epstein \cite{Epstein},
\begin{equation}\label{IntroEpstein}
      \frac{(-1)^{m/2}}{2^{m/2}(2\pi)^{m/2}} \frac{m!}{(m/2)!} \hat{V} = \chi(M), 
\end{equation}
valid for any convex cocompact hyperbolic manifold, which {\em a fortiori} is PE.

In this context, we should mention the recent preprint of Chang, Qing, and Yang \cite{Chang-Qing-Yang}. Their study of Branson's $Q$-curvature allowed them, through a result of Alexakis, to show that
\begin{equation}\label{IntroCQY}
       \int \widetilde{W} \dvol_g 
       + \lrpar{-1}^{\frac{m}{2}} \frac{\Gamma\lrpar{\frac{m+1}{2}}}{\pi^{\frac{m+1}{2}}} \hat{V} 
       = \chi\lrpar{M} .
\end{equation}
In this formula, $\widetilde{W}$ is  a full contraction of the Weyl tensor and its covariant derivatives for the metric $e^{2v}g$ ($e^{v}$ a well-chosen special bdf), and $\hat{V}$ is the renormalized volume. 
Like \eqref{IntroPff} this generalizes Anderson's four-dimensional Gauss Bonnet theorem, as well as Epstein's formula \eqref{IntroEpstein}.
The integrability issues in \cite{Chang-Qing-Yang}
are dealt with by passing to $e^{2v}g$; $(M, e^{2v}g)$ is compact, so the integral of $\widetilde{W}$ needs no renormalization.
A result of Fefferman and Graham, and the choice of $v$, produce the renormalized volume as a boundary term. 

In section $\mathcal{x}$\ref{sec:PEInv}, after rewriting the Pfaffian of an Einstein manifold in terms of the Weyl and scalar curvatures, we note that our formula \eqref{IntroPff} is similarly an integral of a complete contraction of the Weyl curvature plus a multiple of the renormalized volume (the same multiple as in \eqref{IntroCQY}).
In contrast to \eqref{IntroCQY}, we can explicitly identify the Weyl curvature integrand but its integral requires renormalization. The proof of Theorems \ref{Intro1} and \ref{Intro2} depend only on the Fefferman-Graham expansion of the metric (see section $\mathcal{x}$\ref{sec:PE}
 below), whereas that of \eqref{IntroCQY} requires the Einstein condition to apply Alexakis' result.
It would be very interesting to better understand how these formulas are related.

Elsewhere \cite{My thesis}, we consider the heat kernel on a PE manifold and apply renormalization to study index theory. Theorem \ref{Intro1} shows that the coefficients of the small-time asymptotics of the trace renormalize independently of the choice of special bdf. In fact, one can show that the trace of the heat kernel itself, for any fixed positive time, renormalizes in this way. 
Furthermore, using renormalization as $t \to \infty$, we are able to prove an index theorem for the de Rham operator.
It is well-known (see \cite{Mazzeo:Hodge}) that the spaces of harmonic forms on an even-dimensional conformally compact manifold are finite-dimensional except for middle degree forms which always form an infinite-dimensional space. We can use renormalization to define the renormalized index of the de Rham operator ({\em a priori} a real number) and show that:
\begin{equation}\label{HardGB} 
      \Pff = \Rind\lrpar{ L^2\Omega^{\mathrm{even}}(M)%
                 \xrightarrow{d + \delta} L^2\Omega^{\mathrm{odd}}(M) } . 
\end{equation} 
Thus for PE manifolds, \eqref{IntroPff} shows that the renormalized index of the de Rham operator is the Euler characteristic.
The proof  of these results requires much more analysis than that of Theorems \ref{Intro1} and \ref{Intro2} and will not be presented here.

After considering the characteristic numbers of a PE manifold in section $\mathcal{x}$\ref{sec:PEInv}, we study their behavior under variations of the PE structure in section $\mathcal{x}$\ref{sec:Var}.
A family of metrics $g_s$ with Fefferman-Graham expansions (see section $\mathcal{x}$\ref{sec:PE}
 below) imposes a similar expansion on $h$, the infinitesimal variation of $g$, as we show in Proposition \ref{hStructure}.
We analyze the effect of such a variation on the Pontrjagin classes, the renormalized volume and the Pfaffian. For the renormalized volume we recover theorems of Anderson \eqref{intAnd} and Graham-Hirachi \eqref{intGra}. 
We show in Proposition \ref{thm:varPff}, making strong use of the recent work of Labbi \cite{Labbi} on Weyl volume of tubes invariants, that the variation of the renormalized integral of the Pfaffian vanishes. The absence of boundary terms is due to the expansion of $h$, while the vanishing of the interior terms follows from a generalized Bach-Lanczos identity.
Naturally, the vanishing of the variation of the Pfaffian is consistent with equation \eqref{IntroPff}.

\begin{theorem}[\cite{Anderson}, \cite{Graham-Hirachi}] 
A variation of the PE structure on $M$, $g_s$, that preserves the value of the scalar curvature induces a variation of the Pontrjagin characteristic numbers equal to the integral of the appropriate 
Chern-Simons form on the boundary.

The induced variation of the renormalized volume on an even-dimensional manifold is given by
\begin{equation}\label{intAnd}
      \hat{V}'(h) = -\frac{1}{4} \int_{\partial M} \left< g^{(m-1)}, h^{(0)} \right> ,
\end{equation}
while that of the residue of the volume on an odd-dimensional PE manifold is
\begin{equation}\label{intGra} 
     L'(h) = -\frac{1}{4} \int_{\partial M} \left< g^{(m-1,1)},h^{(0)} \right> .
\end{equation}
In these formulas, $g^{(m-1)}$ and $g^{(m-1,1)}$ are the first odd term and first log term, respectively, in the expansion of the metric $g_0$, while $h$ is the infinitesimal variation of $g$ at $s=0$.
\end{theorem} 

The proof of Theorem \ref{Intro2} readily generalizes to other asymptotically regular geometries. For example, in the context of asymptotically Euclidean metrics we obtain a formula for the defect of the topological Euler characteristic and the $L^2$-Euler characteristic involving the Weyl volume of tube invariants of the boundary at infinity. For this and other instances, as well as the proof of \eqref{HardGB}, we refer the reader to \cite{My thesis}.

\subsection{Acknowledgements} 
This work forms part of my thesis. I am very grateful to my advisor, Rafe Mazzeo, for sharing his great erudition and insight. Throughout this work, I received support from his NSF grant DMS-0204730.
I would like to thank Tom Branson for his invitation to participate in the conformal geometry program at the Erwin Schr\"odinger Institute in Vienna in spring of 2004. 
I was fortunate to coincide there with him and Robin Graham, and I would like to thank them both for very interesting conversations.
I am also grateful to the ETH in Zurich for its hospitality during the summer semester of 2004,
and to the Starbucks branches in Zurich and Palo Alto where most of this work was carried out.

\section{Renormalization}\label{sec:MelNis}

\subsection{The Fefferman-Graham expansion} \label{sec:PE}

Suppose that $(M,g)$ is a conformally compact manifold, and $x$ a bdf. It is too optimistic to expect $\bar{g}= x^2g$ to be a smooth metric on $\bar{M}$, and indeed it is well known that $log$ terms arise naturally in its expansion \cite{Fefferman-Graham}. Functions with an expansion of the form
\[ \sum_{k \geq k_0} \sum_{p=0}^{p_k} a_{k,p} x^k \log^p x ,\]
with $a_{k,p}$ smooth functions independent of $x$, are known as ``polyhomogeneous conormal" or phg, see $\mathcal{x}$2A in \cite{Mazzeo:Edge}. We will always assume at least that the metric $\bar{g}$ is continuous and phg. For a discussion of the regularity of the metric, see the recent survey \cite{Anderson3}.

Such metrics are always complete, and any non-trapped geodesic approaches a point on the boundary. 
The sectional curvatures along any geodesic approaching the boundary all converge to $-(|dx|_{\bar{g}})^2$ , and, as a function on the boundary, this is independent of the choice of bdf.
A metric is called asymptotically hyperbolic if $|dx|=1$ on the boundary. These metrics include, and are asymptotically modeled by, hyperbolic metrics. They were introduced in \cite{Mazzeo-Melrose} and \cite{Mazzeo:Hodge} where their resolvents and spectra were studied, respectively.

A choice of bdf determines a metric on the boundary with different choices yielding different but conformally equivalent metrics. Conversely, a boundary metric does not in itself determine a bdf. On an asymptotically hyperbolic manifold, boundary metrics within the conformal class are in one-to-one correspondence with ``special" or ``geodesic" boundary defining functions. A bdf $x$ is special if $|dx|_{\bar{g}}^2 =1$ on a neighborhood of the boundary. The details can be found in Lemma 2.1 of \cite{Graham}; see also \cite{Graham-Lee}.

A choice of special bdf induces, through the flow generated by $\nabla_{\bar{g}}x$, an identification of a neighborhood of $\partial M$ with $\partial M \times [0,\epsilon)$ with metric
\[ \frac{dx^2 + G_x}{x^2}, \]
where $G_x$ is a family of metrics on $\partial M$. Fefferman and Graham \cite{Fefferman-Graham} showed that if $x$ is a special bdf on an $m$-dimensional Poincar\'e-Einstein manifold and $g$ is sufficiently regular, then the expansion of $G_x$ below $x^{m-1}$ is determined by the Einstein condition. Thus if $m$ is even, $\bar{g}$ has an expansion of the form
\begin{equation}\label{Eveng}
         \bar{g} = dx^2 + \bar{g}^{(0)} + x^2\bar{g}^{(2)} + \ldots (\text{even powers}) \ldots 
                                   + x^{m-1}\bar{g}^{(m-1)} + \ldots,
\end{equation}
where $\bar{g}^{(2)},\ldots , \bar{g}^{(m-2)}$ are locally determined by $\bar{g}^{(0)}$ and $\text{tr}_{\bar{g}^{(0)}}\lrpar{\bar{g}^{(m-1)}} = 0$.
For $m$ odd, the analogous expansion is 
\begin{equation}\label{Oddg}
    \bar{g} = dx^2 + \bar{g}^{(0)} + x^2\bar{g}^{(2)} + \ldots (\text{even powers}) \ldots 
    + x^{m-1}\bar{g}^{(m-1)} + x^{m-1}\lrpar{\log x}\bar{g}^{(m-1,1)} + \ldots,
\end{equation}
where now $\bar{g}^{(2)},\ldots , \bar{g}^{(m-3)}$, $\bar{g}^{(m-1,1)}$, and  $\text{tr}_{\bar{g}^{(0)}} \lrpar{ \bar{g}^{(m-1)} }$ are locally determined by $\bar{g}^{(0)}$ and furthermore $\text{tr}_{\bar{g}^{(0)}}\lrpar{\bar{g}^{(m-1,1)}} = 0$.

Following \cite{Graham}, we write the volume form as
\begin{equation}\label{VolForm}
    \text{dvol}_g = \lrpar{ \frac{\det G_x}{\det \bar{g}^{(0)}} }^{1/2} \frac{ \text{dvol}_{\bar{g}^{(0)}} dx}{x^m},
\end{equation}
and note that the expansions above imply
\begin{equation}\label{VolForm2}
     \lrpar{ \frac{\det G_x}{\det \bar{g}^{(0)}} }^{1/2} 
     = 1 + v^{(2)}x^2 + \ldots (\text{even powers}) \ldots + v^{(m-1)}x^{(m-1)} + \ldots,
\end{equation}
where, for $\ell \leq m-1$, $v^{(\ell)}$ is a locally determined function on $\partial M$ and $v^{(m-1)} = 0$ if $m$ is even.

We will say that a phg expansion is $even \mod x^{k}$ if there are no log terms or terms with odd exponents below $x^{k}$.  Thus the metric on a Poincar\'e-Einstein manifold is, in suitable coordinates, $even \mod x^{m-1}$.
Note that the product of two such expansions is again $even \mod x^{k}$.
Graham observed that:

\begin{theorem} \label{thm:Graham}
 If $x$ is a special bdf on a Poincar\'e-Einstein manifold, then 
       \[ \lrpar{ \frac{\det G_x}{\det \bar{g}^{(0)}} }^{1/2} \]
is even$\mod x^{m}$, and the coefficients below  $x^{m}$ are locally determined.
\end{theorem}

We will extend this in Theorem \ref{thm:PEinvs} to full contractions of the curvature and its covariant derivatives. In the rest of this section, we study its implications for renormalization.

As we mentioned above, physical considerations suggested the possibility of renormalizing the volume. This is accomplished by making use of the expansion of the volume form.
Indeed, in even dimensions
\[ \int_{x>\varepsilon} \dvol_g =
 C_0 \varepsilon^{1-m} + C_2 \varepsilon^{3-m} + \ldots (\text{odd powers}) \ldots
 +C_{m-2} \varepsilon^{-1} + \widehat{V} + o(1), \]
and in odd dimensions 
\[ \int_{x>\varepsilon} \dvol_g =
 C_0 \varepsilon^{1-m} + C_2 \varepsilon^{3-m} + \ldots (\text{even powers}) \ldots
 +C_{m-3} \varepsilon^{-2} + L \log \frac{1}{\varepsilon} + \widehat{V} + o(1), \]
where the coefficients $C_{i}$ and $L$ are integrals over $\partial M$ of local curvature expressions of the metric $g^{(0)}$. $\widehat{V}$ is known as the renormalized volume, it depends {\em a priori} on the choice of special bdf.

\subsection{Renormalization Schemes}

The renormalization carried out above, 
\[ \Hint \mu := 
\FP_{\varepsilon =  0} \int_{x > \varepsilon} \mu, \]
 is known as Hadamard regularization. It coincides with the ``b-integral" used in Melrose's proof of the Atiyah-Patodi-Singer index theorem \cite{APS Book}.
 
In a subsequent work, Melrose and Nistor use an alternate method known as Riesz regularization to renormalize integrals (\cite{Melrose-Nistor}; see also \cite{Paycha}). Given a bdf, Riesz regularization of the integral of a density $\mu$ is defined by meromorphically extending 
\begin{equation}\label{zeta}
    \zeta_{x}(z) := \int x^{z}\mu ,
\end{equation}
and taking the finite part at $z=0$; symbolically:
\[ \Rint \mu := \FP_{z=0} \zeta_{x}(z). \] 

Both of these approaches make heavy use of the expansion of the integrand; they are only defined on phg densities. As a preliminary step we assume that the volume form has been written as in \eqref{VolForm} and the boundary integral has been carried out. Thus we only need to consider one-dimensional integrals. Furthermore, making use of linearity, we can localize to a neighborhood of zero, say $[0,\delta)$.

Thus localized, it is easy to compare Hadamard and Riesz renormalizations directly on phg densities. Once the integral along the boundary has been carried out,  
\[ \Hintd x^k \log^p x \text{ dx}= \Rintd x^k \log^p x \text{ dx}
             = \delta^{k+1} \sum_{\ell =0}^p c_{\ell} \log^{p-\ell} \delta, 
             \text{ if $k \neq -1$},\]
while
\[ \Hintd \frac{\log^p x}{x} \text{ dx}= \frac{\log^{p+1} \delta}{p+1} \text{ and }
          \Rintd \frac{\log^p x}{x} \text{ dx}= 0. \]
We can trace this difference back to the fundamental theorem of calculus. Assume that
\begin{equation}\label{fphg}
     f(x) =  \sum_{k \geq k_0} \sum_{p=0}^{p_k} a_{k,p} x^k \log^p x .
\end{equation}
Then
\[ \Hintd f'(x) \text{ dx} =  f(\delta) - a_{0,0}, \text{ and } 
\Rintd f'(x) \text{ dx}= f(\delta) - \sum_{p=0}^{p_0} a_{0,p} \log^p \delta .\]

We define a renormalization scheme $T$ to be a linear functional on phg densities such that, for $f$ as in \eqref{fphg},
\[ T(f'(x) \text{ dx}) = f(\delta) - T^0(f) ,\]
where $T^0(f)$ is a linear function of $\{a_{0,p} \log^p \delta \}$, hence depends on the choice of bdf. 
Given any two renormalization schemes, $T$ and $\widetilde{T}$, it is easy to see that
\begin{equation} \label{TwoTs} \begin{split}
    T\lrpar{ x^k \log^p x \text{ dx}} &= \widetilde{T}\lrpar{ x^k \log^p x \text{ dx}} \text{ for $k \neq -1$, and} \\
    T\lrpar{ \frac{\log^p x}{x} \text{ dx}} - \widetilde{T} \lrpar{ \frac{\log^p x}{x} \text{ dx}} 
    &= \frac{ \widetilde{T}^0 \lrpar{ \log^{p+1} x} - T^0 \lrpar{ \log^{p+1} x} }{p+1}.
\end{split}\end{equation}

The difficulties in comparing 
the same renormalization scheme for two different bdfs
come from the transformation of the coefficients $\{ a_{-1,p} \}$. The situation is greatly simplified by assuming that there are no singular $log$ terms in the expansion of $f$. Note that the absence of singular $log$ terms is independent of the choice of smooth bdf.

\begin{proposition}\label{prop:Res}
Let $\mu$ be density on $M$, phg with respect to $x$, and with no singular log terms. The coefficient $a_{-1,0}$ of $x^{-1}$ in the expansion of $x \mapsto \int_{\partial M_x} \mu$ is given by 
\[ \Res \mu := \Res_{z=0} \int_M x^z \mu \]
and is independent of the choice of bdf. 

If $T$ and $\widetilde{T}$ are any two renormalization schemes, then 
\[T(\mu) - \widetilde{T}(\mu) = C_{T, \widetilde{T}} \Res \mu \]
is independent of the choice of bdf.
\end{proposition}

\begin{proof}
It is easy to see that $a_{-1,0} = \Res_{z=0} \int_M x^z \mu$. 
To see the independence from the choice of bdf, let $\widehat{x}= e^{\omega (x)} x$ be any other bdf. Let $\zeta_x(z)$, $\zeta_{\hat{x}}(z)$ be the zeta functions as defined in \eqref{zeta}. Note that
\[ \zeta_{\hat{x}}(z) - \zeta_x (z)  = \int \lrpar{\hat{x}^z - x^z} \mu 
 = \int \lrpar{e^{\omega z} -1} x^z \mu \\
=  \lrpar{z \int \frac{e^{\omega z}-1}{z} x^z \mu} 
=: z \widetilde{\zeta}(z) .\]
Due to the absence of singular $log$ terms, the meromorphic continuation of $\widetilde{\zeta}$ has at most a simple pole at $z=0$. Hence $\zeta_{\hat{x}} - \zeta_{x}$ extends to be holomorphic at $z=0$, and the residues must be equal.

The final statement follows directly from \eqref{TwoTs}.
\end{proof}

In what follows, we will only seek to renormalize integrals of densities without singular $log$ terms, and thus we shall unabashedly consider only Riesz renormalization, allowing for much simpler computations. A case in point is the proof of Proposition \ref{prop:Res} which yields the following proposition, taken from \cite{Melrose-Nistor}.

\begin{proposition} \label{prop:2bdfs}
Let $\mu$ be a density without singular log terms, and let $x$ and $\widehat{x}= xe^{\omega(x)}$ be two bdfs.
The difference between the renormalized integrals with respect to $\hat{x}$ and $x$ is given by
\begin{equation}\label{difint}
     \widehat{\Rint \mu} - \Rint \mu = \int_{\partial M} \lrpar{\omega\mu}_{(-1)}
\end{equation}
where $\lrpar{ \omega \mu }_{(-1)}$ denotes to the term in the expansion of $ \omega \mu$ of homogeneity $-1$ in $x$.
\end{proposition}

\begin{proof}
Simply note from the above computation that
\begin{equation}\label{difzetas}
            \FP_{z=0} \lrpar{\zeta_{\hat{x}}(z) - \zeta_x(z)}
      = \FP_{z=0} \lrpar{z \int \frac{e^{\omega z}-1}{z} x^z \mu} 
      = \Res_{z=0} \int x^z\omega \mu
        = \int_{x=0} \lrpar{ \omega \mu }_{(-1)}.
\end{equation}
\end{proof}

\subsection{Renormalization on Poincar\'e-Einstein manifolds}

In general,
for a density to admit a renormalization independent of the choice of bdf, it needs to actually be integrable. Nevertheless, we will see that on Poincar\'e-Einstein manifolds there is a rich class of densities renormalizing independently of the choice of {\em special} bdf. The reason is twofold. On the one hand, the metric has the Fefferman-Graham expansion \eqref{Eveng}, \eqref{Oddg}. On the other hand, we have the following lemma from \cite{Graham-Lee} (see also \cite{Guillarmou:Mero Prop of Resol}).

\begin{lemma}\label{lem:omega}
Let $(M,g)$ be an asymptotically hyperbolic manifold, and $x$ a special bdf. 
For any odd number $k$, the expansion of $x^2g$ is even below $x^k$ if and only if for any other special bdf $\hat{x}=e^{\omega(x)}x$ the expansion of $\omega$ is even below $x^{k+2}$.
\end{lemma}

\begin{proof} 
 The condition $|d\hat{x}|^2=1$ imposes
\begin{equation}\label{dx=1}
          2\partial_x\omega + x\lrpar{(\partial_x\omega)^2
          + \bar{g}^{ij}(\partial_i\omega)(\partial_j\omega)} = 0.
\end{equation}
In terms of the expansion of $\omega$ in $x$, this does not restrict $\omega^{(0)}$. On the other hand, it is easy to see that if $\omega$ is even below $x^r$, then the first term in \eqref{dx=1} is odd below $x^{r-1}$, while the other terms are odd below $x^{min(k+1, r+1)}$, hence for $r<k+2$, $\omega^{(r)}=0$.
Thus the first odd term is $\omega^{(k+2)}$, and it is given by:
\[ (k+2)\omega^{(k+2)} = -\lrpar{\bar{g}^{(k)}}^{ij} (\partial_i\omega^{(0)})(\partial_j\omega^{(0)}).\]
This is zero for arbitrary $\omega^{(0)}$ precisely when $\bar{g}^{(k)} =0$.
\end{proof}

Recall from \eqref{VolForm} and \eqref{VolForm2}, that the volume form of an even-dimensional PE manifold has an expansion of the form
\[ \dvol_g = \lrpar{ 1 + v^{(2)}x^2 + \ldots \text{(even powers)} \ldots + v^{(m)}x^{m} + \ldots }
                  \frac{\dvol_{g^{(0)}} dx}{x^m} .\]
Note that $\dvol_g$ has no residue, and that the lemma guarantees that changing to another special bdf will not produce a residue.
In view of Proposition \ref{prop:2bdfs}, this means that the renormalized integral of the volume is independent of the choice of special bdf. This has the following immediate generalization.

\begin{theorem}\label{thm:ReInt}
Let $(M^m, g)$ be a Poincar\'e-Einstein manifold, $x$ a special bdf, and $\mu$ a density of the form
        \[ \mu = F(x,y) \lrpar{ \frac{\dvol_{g^{(0)}} dx}{x^m} }, \]
where $F$ is $even \mod x^{m+1}$. Then, if $m$ is even, $\Rint \mu$
is independent of the choice of special bdf.
If $m$ is odd, and $\widehat{x}=xe^{\omega(x)}$ is another special bdf, then
       \[ \widehat{\Rint \mu} - \Rint \mu = \int_{\partial M} \lrpar{ \omega F(x,y) }_{(m-1)}. \]
Though in any case, $\Res \mu$ is independent of the choice of bdf.
\end{theorem}

One might initially expect the renormalized volume, coming from a conformal theory, to depend only on the conformal class. Its definition involves the choice of a bdf and this breaks the conformal invariance, explaining the appearance of a conformal anomaly in odd dimensions. In even dimensions the conformal invariance is restored after renormalization.

To apply this theorem to the scalar Riemannian invariants of a PE manifold, we need to understand the expansions of their curvature tensors.

\section{Curvature of an Asymptotically Hyperbolic Manifold}\label{sec:AH}

In this section we develop some of the geometry of an asymptotically hyperbolic manifold. Though our applications in the next section will all be to PE manifolds, we will not assume that the metric is Einstein in this section. We begin by introducing the coordinates in which our computations will be carried out.

Let $(M^m,g)$ be an asymptotically hyperbolic manifold, and $x$ a special bdf. Pick a point on the boundary $p \in \partial M$, and choose coordinate vector fields $\{ \partial_{y_i}, \partial_x \} =: \{ \bar{X}_s \}$ for $\bar{g}$ by exponentiating first on the boundary, then into the manifold.  That is, $\{ \partial_{y_i} \}$ form a normal coordinate chart for $(\partial M,\bar{g})$ centered at $p$ and are extended into the interior of $M$ along geodesics normal to $\partial M$.
In this way, with $x$ as the $m^{\mathrm{th}}$ coordinate, throughout the chart we have
\begin{equation}\label{Orth}
    \bar{g}_{km} = \delta_{km}. 
\end{equation}
The Christoffel symbols and components of the curvature tensor of $\bar{g}$ in this coordinate chart will be denoted by $\bar{\Gamma}^s_{tu}$ and $\bar{R}^s_{tuv}$, respectively. 
We shall consider the structure of $g$ using the frame $X_s := x\bar{X}_s$. We will use the letters $i$, $j$, $k$, and $\ell$ as indices varying between $1$ and $m-1$ and $s$, $t$, $u$, and $v$ to denote indices between $1$ and $m$. Note that $\lrspar{ \bar{X}_s, \bar{X}_t } = 0$ for all $s,t$, but
\begin{equation}\label{Tor}
    \lrspar{ X_m, X_i} = X_i , \lrspar{X_i,X_j}=0, \text{ for }i,j<m.
\end{equation}
Also, in this frame
\begin{equation}\label{ThisFrame}
 g\lrpar{ X_s, X_t } = x^2g\lrpar{ \bar{X}_s, \bar{X}_t} = \bar{g}_{st}
 \end{equation}

The Levi-Civita connection of $g$ is closely related to that of $\bar{g}$.
Indeed, Theorem 1.159 of \cite{Besse} specializes to:
\begin{equation}\label{Besse}
      \nabla_W Y = \bar{\nabla}_W Y
      - \frac{dx}{x}\lrpar{W} Y - \frac{dx}{x}\lrpar{Y}W + \frac{\bar{g}\lrpar{W,Y}}{x^2} x\partial_x .
\end{equation}
In terms of the frame $X_i$ above, in analogy with the Christoffel symbols, we define $\gamma^u_{st}$ by the equation $\nabla_{X_s}X_t = \gamma^u_{st} X_u$.

\begin{lemma}\label{lem:Gammas}
With $X_s = x\bar{X}_s$ as above,
\begin{equation}\label{gammaeq}
       \gamma^t_{us} = x\bar{\Gamma}^t_{us} - \delta_{sm}\delta_{ut} + \delta_{tm}\bar{g}_{us} .
\end{equation}
Furthermore, recalling \eqref{Orth},
\begin{equation*}
       \gamma^m_{ij} = -\frac{1}{2} \lrpar{ x\partial_x \bar{g}_{ij} } + \bar{g}_{ij} \text{ , }
       \gamma^k_{im} = \frac{1}{2} \bar{g}^{k \ell} \lrpar{ x\partial_x \bar{g}_{i\ell} }  - \delta_{ik} \text{ , }
       \gamma^k_{mi} = \frac{1}{2} \bar{g}^{k \ell} \lrpar{ x\partial_x \bar{g}_{\ell i} }.
\end{equation*}
\end{lemma}

Moving on to the curvature, define $r^s_{tuv}$ by the equation
        \[ R(X_u,X_v)X_t = r^s_{tuv} X_s ,\]
and, directly from $R(X,Y)Z = \lrpar{ \lrspar{\nabla_X, \nabla_Y} - \nabla_{\lrspar{X, Y}} }Z$, obtain
\begin{equation}\label{Curv}
        r^s_{tuv} = \gamma^s_{uw}\gamma^w_{vt} - \gamma^s_{vw}\gamma^w_{ut}
                           + X_u \lrpar{\gamma^s_{vt} } - X_v \lrpar{\gamma^s_{ut} }
                           -\delta_{um} \gamma^s_{vt} + \delta_{vm} \gamma^s_{ut}. 
\end{equation}

An expansion of $g$ in $x$ induces expansions of $\gamma^s_{tu}$ and $r^s_{tuv}$. Note that the expansion of $\bar{\Gamma}^s_{tu}$ in $x$ below $x^n$ uses the expansion of $\bar{g}$ below $x^{n+1}$. An advantage of using $x\partial_x$ instead of $\partial_x$ is that it does not lower order of homogeneity and hence, for example, the expansion of $\gamma^s_{tu}$ below any $x^n$ depends only on the expansion of $\bar{g}$ below $x^n$. Indeed, this is still true of the curvature and even its covariant derivatives. We shall prove this with an eye towards Poincar\'e-Einstein metrics.

\begin{theorem}\label{thm:CurvJets}
Let $(M^m,g)$ be an asymptotically hyperbolic manifold, and $x$ a special bdf.
If $\bar{g}$ has a phg expansion that extends continuously to $\bar{M}$ (thus $x^k \log x$ can occur for $k \geq 1$, but not for $k=0$),
then so does the curvature and any of its covariant derivatives. Furthermore, the terms below $x^n$ in the expansion of any covariant derivative of the curvature are determined by the terms below $x^n$ in the expansion of $\bar{g}$.
\end{theorem}

{\em Remark.} A consequence of using a frame for $g$ built from a frame for $\bar{g}$, i.e., $X_s = x\bar{X}_s$ is that we raise and lower indices using $\bar{g}$ instead of $g$. For instance,
\begin{equation*}\begin{split}
 r_{stuv} &= g\lrpar{ R\lrpar{ X_u, X_v }X_t, X_s } 
                  = g\lrpar{ r^w_{tuv} X_w, X_s} 
                  = \lrpar{r^w_{tuv} } g\lrpar{ x\bar{X}_w, x\bar{X}_s} 
                  = \bar{g}_{sw} r^w_{tuv}.
 \end{split}\end{equation*}
Hence the conclusion of the theorem does not change if we raise or lower indices. 

\begin{proof}
As stated before the theorem, this is a simple consequence of using $x\partial_x$ instead of $\partial_x$.
Note that $x\partial_x \lrpar{x^k} = kx^k$ and that 
$x\partial_x \lrpar{x^k \log^p x} = k x^k \log^p x + p x^k \log^{p-1}x$, so the formulas \eqref{gammaeq}, \eqref{Curv}
show that the theorem holds for $\gamma^s_{tu}$ and $r^s_{tuv}$. For the covariant derivatives of the curvature we use induction and the formula
\begin{equation}\label{DR}
 \lrpar{\nabla^{p+1}R}_{stuv; \alpha_1 \cdots \alpha_p, \alpha_{p+1}}
     = X_{\alpha_{p+1}} \lrpar{ r_{stuv; \alpha_1 \cdots \alpha_p}}
        - \gamma^{\beta}_{\alpha_{p+1} s} r_{\beta tuv; \alpha_1 \cdots \alpha_p}
        -\ldots - \gamma^{\beta}_{\alpha_{p+1} \alpha_{p}} r_{stuv; \alpha_1 \cdots \beta}.
\end{equation}

\end{proof}

For the PE context, one is interested in a metric with an expansion
containing only even powers of $x$ below some term, say $x^n$. To organize our discussion of such metrics, define $\mathcal{F}$ of a function to be $1$ if the function is $even$ below $x^{n}$ and $-1$ if it is $odd$ below $x^{n}$. $\mathcal{F}$ is clearly a multiplicative homomorphism among functions with such expansions.

\begin{corollary} \label{cor:Even}
In the above context, assume that the expansion of $\bar{g}$ below $x^n$ consists entirely of even exponents.  Then the expansion of each $r_{stuv; \alpha_1 \cdots \alpha_p} $ is either even or odd below $x^n$. Indeed, if $q$ is the number of $m$'s among the indices of $r_{stuv; \alpha_1 \cdots \alpha_p}$, it satisfies
\begin{equation}\label{FofR}
     \mathcal{F}(r_{stuv; \alpha_1 \cdots \alpha_p}) = (-1)^{p+q}. 
\end{equation}
If $P$ is a full contraction of the curvature and its covariant derivatives (i.e., a scalar Riemannian invariant), then $\mathcal{F}(P) = 1$.
\end{corollary}

\begin{proof}
Because of Theorem \ref{thm:CurvJets}, we may truncate the expansion of $\bar{g}$ below $x^n$ so that it consists entirely of even powers of $x$. 
Our formulas for $\gamma^s_{tu}$ in Lemma \ref{lem:Gammas} show that
\[ \mathcal{F}\lrpar{\gamma^s_{tu}} = 
\begin{cases}
1 & \text{ if $m \in \{s,t,u\}$} \\
-1 & \text{ otherwise}
\end{cases}. \]

For any function $f$,
\[ \mathcal{F}\lrpar{ X_s \lrpar{f} } = 
\begin{cases}
 \mathcal{F}\lrpar{f} & \text{ if $s= m$} \\
 -\mathcal{F}\lrpar{f} & \text{ otherwise}
 \end{cases}, \]
so a close look at the formula for the curvature \eqref{Curv} shows that the parity of $r^s_{tuv}$ is just the parity of the number of $m$'s among $\{s,t,u,v\}$. 
As we have chosen a chart where $\bar{g}_{jm} = \delta_{jm}$,
\[ \mathcal{F}\lrpar{r^s_{tuv}} = \mathcal{F}\lrpar{r_{stuv}}, \]
which proves the theorem in the case $p=0$.

Induction and the formula \eqref{DR} finish the proof of \eqref{FofR}.
Note that once we know \eqref{FofR} it is clear that raising or lowering indices in the curvature or any of its covariant derivatives does not change the value of $\mathcal{F}$. 

Consider now $P$, a full contraction of the curvature and its covariant derivatives.
We know that $\mathcal{F}$ of any such contraction is the product of $\mathcal{F}$ of the individual factors, each of which is $(-1)$ raised to the number of derivatives plus the number of $m$'s among the indices. Hence $\mathcal{F}$ of the full contraction is $(-1)$ raised to the total number of $m$'s appearing among the indices plus the total number of derivatives among all of the factors. The total number of derivatives is necessarily even, as is the total number of $m$'s among the indices, since the indices are paired together to form the contraction. This proves that any full contraction has only even exponents below $x^{n}$.
\end{proof}

\section{Invariants of a Poincar\'e-Einstein Manifold}\label{sec:PEInv}

\subsection{Scalar Riemannian Invariants}\label{ScalInv}
Throughout this section $(M,g)$ is a Poincar\'e-Einstein manifold, and $x$ will always denote a special bdf. Recall from section $\mathcal{x}$\ref{sec:PE} that, in this context, the expansion of $\bar{g}=x^2g$ is even below $x^{m-1}$ (regardless of the parity of $m$). Let $P$ be any scalar Riemannian invariant (i.e., a full contraction of the curvature and its covariant derivatives). Corollary \ref{cor:Even} guarantees that $P$ is also even below $x^{m-1}$. Our first task is to understand the effect of the $x^{m-1}$ term in the expansion of $\bar{g}$ on the $x^{m-1}$ term in expansion of $P$. 

Recall that the frame $\{X_u\}$ is centered at a point $p \in \partial M$. We may assume that, at $p$, $g^{(0)}$ is the identity and $g^{(m-1)}$ is diagonal, say $g^{(m-1)}_{ij} = \mu_i \delta_{ij}$. We will work under these assumptions and compute the expansion of $P$ at the point $p$.

\begin{lemma}\label{lem:m-1}
Let $(M^m,g)$ be an asymptotically hyperbolic manifold, and $x$ a special bdf.
Assume the metric has an expansion of the form
\[    \bar{g} = dx^2 + \bar{g}^{(0)} + x\bar{g}^{(1)} + \ldots + x^{n-1}\bar{g}^{n-1} 
    + x^{n}\bar{g}^{(n)} + x^{n}\lrpar{\log x}\bar{g}^{(n,1)} + \ldots, \]

and let $P$ be a full contraction of the curvature and its covariant derivatives. Then It has a similar expansion
      \[ P = P^{(0)} + xP^{(1)} + \ldots + x^{n-1}P^{(n-1)}
    + x^{n}P^{(n)} + x^{n}\lrpar{\log x}P^{(n,1)} + \ldots, \]
wherein the dependence of $P^{(n)}$ and $P^{(n,1)}$ on $\bar{g}^{(n)}$ and $\bar{g}^{(n,1)}$ is only through
$trace_{\bar{g}^{(0)}}\bar{g}^{(n)}$ and $trace_{\bar{g}^{(0)}}\bar{g}^{(n,1)}$.
\end{lemma}

\begin{proof}
For our purpose, we may consider $\bar{g}$ given by:
\[ dx^2 + \bar{g}^{(0)} + \bar{g}^{(n)} x^n + \bar{g}^{(n,1)} x^n \log x .\]
Indeed, note that we do not need to compute the $x^{n}$ term in the expansion of $P$, only its contribution from the $x^n$ term in $\bar{g}$, and interactions between say $\bar{g}
^{(n)}$ and $\bar{g}^{(i)}$ do not occur below $x^{n+i}$.

We will compute at the point $p \in \partial M$ as discussed before the statement of the lemma. Say that a function, $f$, is {\em respectable} if the $x^0$ term is a constant and the contribution of the $x^n$ term in the expansion of $\bar{g}$ to the $x^n$ term in the expansion of $f$ is only through a (constant) linear combination of the components of $\bar{g}^{(n)}$ and $\bar{g}^{(n,1)}$. It is easy to see from the formulas in Lemma \ref{lem:Gammas} that $\gamma^s_{tu}$ is respectable.

Note that if $f$ is respectable, then so are $x\partial_y(f)$ and $x\partial_x(f)$. Also, the product of respectable functions is again respectable.
It follows that any component of the curvature computed in such a frame is respectable, and similarly, components of covariant derivatives of the curvature are respectable.

Thus in the contraction we are interested in, the individual terms are respectable. We finish by noting that contracting a linear combination of components of $\bar{g}^{(n)}$ and $\bar{g}^{(n,1)}$ yields a linear combination of $trace_{\bar{g}^{(0)}}\bar{g}^{(n)}$ and $trace_{\bar{g}^{(0)}}\bar{g}^{(n,1)}$.
\end{proof}

We can now state our generalization of Theorem \ref{thm:Graham}.

\begin{theorem}\label{thm:PEinvs}
If $x$ is a special bdf on a PE manifold, then any full contraction of the curvature and its covariant derivatives has an even expansion$\mod x^{m}$, and the coefficients below $x^{m}$ are locally determined.

The behavior of scalar Riemannian invariants under renormalization is thus given by Theorem \ref{thm:ReInt}. In particular, on an even-dimensional PE manifold, all scalar invariants renormalize independently of the choice of special bdf.
\end{theorem}

\begin{proof}
Recalling the expansions \eqref{Eveng} and \eqref{Oddg},
Theorem \ref{thm:CurvJets} establishes the theorem below $x^{m-2}$, and we use Lemma \ref{lem:m-1} to deal with the $x^{m-1}$ term.
\end{proof}

As is clear from the proof of the theorem, the Einstein condition only enters through the Fefferman-Graham expansion and the result holds whenever $x^2g$ is even below $x^{m-1}$ as long as \newline
 $trace_{g^{(0)}}g^{(m-1)} =0$.

Also recall, from Proposition \ref{prop:2bdfs}, that changing from one special bdf $x$ to another one, $\hat{x}=e^{\omega}x$, changes the volume by $\int (\omega \dvol)_{-1}$. For a PE metric, this term vanishes due to the Fefferman-Graham expansion of the metric, but if the metric were ``any less even" this integrand would certainly not be zero in general.
Even expansions of asymptotically hyperbolic (not necessarily Einstein) metrics are related in \cite{Guillarmou:Mero Prop of Resol} to the domain on which the resolvent is meromorphic with poles of finite rank. It would be interesting to find a direct connection between the extent of this domain and the renormalization of the volume, for example, without going through the expansion of the metric.

\subsection{Pontrjagin Classes}\label{Pont}

We start by reviewing the conformal invariance of the Pontrjagin classes. In dealing with these, it is useful to think of the curvature as an endomorphism-valued 2-form. On any manifold, define the 2-form ${\mathbf \Omega}^u_v$ by
\[ {\mathbf \Omega}^u_v = R^u_{vst} \theta^s \wedge \theta^t ,\]
where $R$ is computed with respect to a frame $\{ X_i \}$ with dual frame $\{ \theta^i \}$.
In the same way, we define the 2-form $\mathbf{W}^u_v$ replacing the components of $R$ above with the corresponding components of the Weyl curvature tensor. It is easy to see that the forms $\mathbf{W}^u_v$ are conformally invariant 2-forms.

The following theorem was first shown by Avez \cite{Avez}, though the conformal invariance of the Pontrjagin classes had already been established by Chern and Simons.
\begin{theorem}\label{Avez}
On any Riemannian manifold
\[ \tr\lrpar{ {\mathbf \Omega}^{\ell} } = \tr\lrpar{ {\mathbf W}^{\ell}} \]
\end{theorem}

In particular, on a conformally compact manifold, any characteristic number built up from Pontrjagin classes is integrable and equals the corresponding characteristic number of $(\bar{M},\bar{g})$.

\subsection{The Pfaffian}\label{sec:Pff}

On an even-dimensional manifold, we have one more characteristic class, the Pfaffian. As a PE scalar Riemannian invariant, we know that its integral renormalizes independently of the choice of special bdf, and in fact (Theorem \ref{SoftGB} below)
\begin{equation} \Pff = \chi(M) . \end{equation}

As with any scalar Riemannian invariant on an Einstein manifold, the Pfaffian is a polynomial in the scalar and Weyl curvatures. We will identify this polynomial in terms of the ``volume of tubes" invariants of the Weyl curvature and begin by reviewing these.

Weyl (see \cite{Tubes}) derived a formula for the volume of a tube of small radius $\varepsilon$ around a $q$-dimensional submanifold $P$ of $\mathbb{R}^n$:
\[ Vol_P(\varepsilon) = \omega_{n-q}  \sum_{\ell=0}^{\lfloor q/2 \rfloor} 
        \frac{K_{2\ell}(P) \varepsilon^{2\ell}}{(n-q +2)(n-q+4)\cdots(n-q+2\ell)} ,\]
 where $\omega_{n-q}$ is the volume of the unit ball in $\mathbb{R}^{n-q}$. 
Weyl showed that the coefficients are intrinsic to $P$ by identifying them with integrals over $P$ of complete contractions of the curvature and its covariant derivatives.  Donnelly \cite{Donnelly} later proved that they are universal linear combinations of the heat invariants of $P$.
 
Weyl's theorem was used in the first proofs of the Gauss Bonnet theorem for a general even-dimensional compact manifold (independently by Allendoerfer and Fenchel for submanifolds of $\mathbb{R}^n$, then by Allendoerfer and Weil in general). The link is 
 \[K_{\dim P}(P) = (2\pi)^{(\dim P)/2}\chi(P) .\]

One expression for $K_{2\ell}$ comes from considering the curvature as an endomorphism valued 2-form as in $\mathcal{x}$ \ref{Pont}. Indeed, the coefficients are multiples of the Lipschitz-Killing curvatures (see \cite{Cheeger})
\begin{equation}\label{LipKil} 
K_{2\ell} = \int C_{\ell} \sum_{\sigma \in S_n} \text{sign}(\sigma) \phantom{x}
   \Omega_{\sigma_1}^{\sigma_2} \wedge \cdots \wedge  \Omega_{\sigma_{2\ell-1}}^{\sigma_{2\ell}} 
   \wedge \theta^{\sigma_{2\ell+1}} \wedge \cdots \wedge \theta^{\sigma_m} .
\end{equation}

To obtain another expression for $K_{2\ell}$, it is convenient to work in the formalism of double forms as set out in \cite{Tubes} and especially \cite{Kulkarni}. In this way we will obtain the aforementioned polynomial for the Pfaffian.

Define
\[ \curly{D}^{p,q} := \Lambda^p(M) \otimes \Lambda^q(M). \] 
We consider the metric and the curvature as (symmetric) elements of $\curly{D}^{1,1}$, $\curly{D}^{2,2}$, respectively. They are given by
\[ \df{g}(X \otimes Y) = g(X,Y), \text{ and }  \df{R} \lrpar{(X,Y) \otimes (Z,W)} = g(R(X,Y)Z,W). \]

There is a natural operation on $\curly{D}^{*,*}$, the Kulkarni-Nomizu product given by
\[ \lrpar{\df{a}_1 \otimes \df{a}_2} \cdot \lrpar{\df{b}_1 \otimes \df{b}_2}
                           = \df{a}_1 \wedge \df{b}_1 \otimes \df{a}_2 \wedge \df{b}_2 .\]
For instance the curvature tensor with constant sectional curvature, $\lambda$, is given by
$ \frac{\lambda}{2} \df{g}^2$. Note that this is sometimes, e.g. Definition 1.110 of \cite{Besse}, denoted $\owedge$.
                           
We also introduce the contraction map $\curly{C}: \curly{D}^{p+1,q+1} \to \curly{D}^{p,q}$. Let $E_i$ be any locally defined orthonormal frame; then for $\df{P} \in \curly{D}^{p+1,q+1}$,
\[ \curly{C}\df{P} \lrpar{(X_1,\ldots,X_p) \otimes (Y_1,\ldots,Y_q)} :=
               \sum_{\ell=1}^{m} \df{P}\lrpar{(E_{\ell},X_1,\ldots,X_p) \otimes (E_{\ell},Y_1,...,Y_q)} .\]
Thus $\curly{C}\df{R}$ is the Ricci curvature and $\curly{C}^2\df{R}$ is the scalar curvature.

The volume of tube invariants are given by
\begin{equation}\label{VolInv}
          K_{2\ell} = \int \frac{\curly{C}^{2\ell}(\df{R}^{\ell})}{(2\ell)!\ell !}  \dvol .
\end{equation}
Lemma 5.5 in \cite{Tubes} expresses the Pfaffian of an $m$-dimensional manifold as
\begin{equation}\label{Pff1}
     \text{Pff} = \frac{1}{(2\pi)^{m/2}} \frac{ \curly{C}^{m}\df{R}^{m/2}}{m!(m/2)!} \dvol.
\end{equation}

The Weyl curvature is defined by
\[ \df{W}= \df{R}- \frac{ \df{g}  \cdot \curly{C}\df{R} }{m-2} 
              + \frac{ \df{g}^2 \cdot \curly{C}^2\df{R}}{2(m-1)(m-2)} ,\]
so on an Einstein manifold
\[ \df{R} 
              = \df{W} + \frac{s_g}{2m(m-1)} \df{g}^2 . \]

Note that since $\df{W}$ and $\df{g}^2$ are in $\curly{D}^{2,2}$, a commutative algebra, the binomial theorem applies to give
\begin{equation}\label{Pff2}
 \df{R}^{m/2} =
      \sum_{k=0}^{m/2} \binom{m/2}{k} \lrpar{\frac{s_g}{2m(m-1)}\df{g}^2}^k \cdot \df{W}^{\frac{m}{2}-k} .
\end{equation}
Furthermore, for any $\df{P} \in \curly{D}^{\ell,\ell}$ we have the following formula from \cite{Kulkarni}:
\[ \curly{C}^{\ell+1} \lrpar{\df{g} \cdot \df{P}} 
        = (\ell +1)(m-\ell)\curly{C}^{\ell}\df{P} .\]
Iterating we get
\begin{equation}\label{Pff3}
 \curly{C}^{m} \lrpar{\df{g}^{2k} \cdot \df{W}^{\frac{m}{2}-k}} 
        =  \frac{m!(2k)!}{(m-2k)!} \curly{C}^{m-2k} \df{W}^{\frac{m}{2}-k} .
\end{equation}

Finally, plugging \eqref{Pff2} and \eqref{Pff3} back into \eqref{Pff1} we obtain the following formula.

\begin{lemma}\label{EinPfaff}
On an $m$ dimensional Einstein manifold, the Pfaffian is given by
\begin{equation*}
           \frac{1}{(2\pi)^{m/2} }
           \sum_{k=0}^{m/2} \frac{(2k)!}{k!}  
           \lrpar{\frac{s_g}{2m(m-1)}}^k
           \frac{\curly{C}^{m-2k} \df{W}^{\frac{m}{2}-k} } { \lrpar{m-2k}! \lrpar{\frac{m}{2} -k}!} \dvol
\end{equation*}
\end{lemma}

Note that the factors $\curly{C}^{m-2k} \df{W}^{\frac{m}{2}-k} $ are (multiples of) the Weyl volume of tubes invariants, except that they are evaluated in the Weyl curvature instead of the full curvature tensor.
Note the scalar curvature factor in each summand with $k \neq 0$; these terms are not conformally invariant.

Returning to our PE manifold, we may replace $s_g = -m(m-1)$, and
this lemma gives us a formula for the Pfaffian
 in terms of the successive contractions of the Weyl tensor. 
Although the Weyl tensor itself has trivial contraction, this is not true of its higher powers. For instance, from \cite{Tubes}, $\curly{C}^4 \df{W}^2 = 6|W|^2 $. 

\begin{theorem}\label{SoftGB}
Let $x$ be a special bdf on an even dimensional asymptotically hyperbolic manifold $(M,g)$, and assume that $x^2g = dx^2 + G_x$ where
     \[ G_x = \bar{g}^{(0)} + ...(\text{even powers})... + \bar{g}^{(m-1)} x^{m-1} + O(x^m) ,\]
and $trace_{\bar{g}^{(0)}} \bar{g}^{(m-1)} =0$. Then
     \[ \Pff = \chi(M) .\]
\end{theorem}

\begin{proof}
This follows from Chern's Gauss-Bonnet Theorem for the manifold with boundary
$M_{\varepsilon} := \{x \geq \varepsilon \}$,
\begin{equation}\label{ChernGB}
      \int_{M_{\varepsilon}} \text{Pff} + \int_{x=\varepsilon} \Fun = \chi(M_{\varepsilon}). 
\end{equation}
Here $\Fun$ is a polynomial in the curvature and the second fundamental form of $\{x = \varepsilon \}$ in 
$M_{\varepsilon}$. In terms of the curvature and connection forms it is given in \cite{Chern} by
\begin{equation}\label{ChernFun}
      \Fun = \sum_{q=0}^{m/2} C_{m,q} 
       \sum_{\sigma \in \Sigma_{m-1}} \Omega_{\sigma_1\sigma_2} \wedge \ldots
            \wedge \Omega_{\sigma_{2q-1}\sigma_{2q}} \wedge \omega_{\sigma_{2q+1}m}
            \wedge \ldots \wedge \omega_{\sigma_{m-1}m} ,
\end{equation}
for some constants $C_{m,q}$ whose precise value we will not need. Note that, for $\varepsilon$ small, $\chi(M_{\epsilon}) = \chi(M)$. Since the right hand side of \eqref{ChernGB} does not depend on $\varepsilon$, neither does the left, and thus
\[ \Pff + \FP_{\varepsilon = 0} \int_{x=\varepsilon} \Fun = \chi(M). \]
So we need only show that the second term vanishes.

If we denote the second fundamental form of $\{x = \varepsilon \}$ as a double form by $\FunForm \in \curly{D}^{1,1}$, \eqref{ChernFun} shows that $\Fun$ is a linear combination of terms of the form
\begin{equation}\label{FunTerm}
      \curly{C}^{m-1} (\df{R}^k\FunForm^{m-1-2k}) \dvol.
\end{equation}
It is easy to see that, for any $\omega$,
\[ \curly{C}^{m-1}(\omega) \dvol_g\downharpoonright_{\partial M} 
        = x^{m-1} \bar{\curly{C}}^{m-1}(\omega) \dvol_{\bar{g}}\downharpoonright_{\partial M} ,\]
where $\bar{\curly{C}}$ denotes contracting via a local orthonormal frame of $\bar{g}$ instead of $g$.
Recall that $r_{ijk\ell}$ and $\gamma^m_{ij}$ have even expansions mod $x^{m-1}$, hence so do $x^4\df{R}$ and $x^2\FunForm$, and the constant term in \eqref{FunTerm} comes from the $x^{m-1}$ terms in the expansion of these coefficients. Now, since
\begin{equation*}\begin{split}
      \lrpar{\gamma^m_{ij}}^{(0)}= \bar{g}_{ij}^{(0)}  &, %
        \lrpar{\gamma^m_{ij}}^{(m-1)}= \frac{3-m}{2}\bar{g}_{ij}^{(m-1)} , \\
         \lrpar{r_{ijk\ell}}^{(0)} =%
              \lrspar{ \gamma^m_{ik}\gamma^m_{j\ell} - \gamma^m_{jk}\gamma^m_{i\ell} }^{(0)} &,%
        \text{ and }%
         \lrpar{r_{ijk\ell}}^{(m-1)} =%
               \lrspar{ \gamma^m_{ik}\gamma^m_{j\ell} - \gamma^m_{jk}\gamma^m_{i\ell} }^{(m-1)} ,
\end{split}
\end{equation*}
we may conclude that the constant term in the expansion of \eqref{FunTerm} is a multiple of \newline $\bar{\curly{C}}^{m-1}(\bar{\df{g}}_{(0)}^{m-2}\bar{\df{g}}_{(m-1)}) \bar{\dvol}$, which vanishes as 
$trace_{\bar{g}^{(0)}} \bar{g}^{(m-1)} =0$.
\end{proof}

It is instructive to use Lemma \ref{EinPfaff} and write out the renormalized Gauss-Bonnet theorem \ref{SoftGB} in four dimensions:
\begin{equation}\label{Anderson}
  \frac{1}{(2\pi)^{2}} \Rint
           \sum_{k=0}^{2} \frac{(2k)!}{k!}
           \lrpar{-\frac{1}{2}}^k \frac{\curly{C}^{4-2k} \df{W}^{2-k}}{(4-2k)!(2-k)!}
           = \frac{1}{8(2\pi)^{2}} \int |W|^2
           + \frac{3}{(2\pi)^{2}}  \hat{V}  = \chi(M) ,
\end{equation}
and in six dimensions:
\begin{equation}\label{6d}
      \frac{1}{(2\pi)^{3}} \int \frac{\curly{C}^{6} \df{W}^{3}}{6!3!}
               -\frac{1}{8(2\pi)^{3}} \Rint |W|^2
               -\frac{15}{(2\pi)^{3}} \hat{V} = \chi(M).
\end{equation}
Equation \eqref{Anderson} is due to Anderson \cite{Anderson} 
(note that he uses a different convention for $|W|^2$, see the comment after equation (1.15) in \cite{Anderson}).
It is only in dimension four that the integrand requires no further renormalization than  the volume. Generally, from Lemma \ref{EinPfaff} we know that
\begin{equation}\label{tocomp}
 \Pff = \int \text{Pff}(W) + \Rint P(W) + 
           \frac{(-1)^{m/2}}{2^{m/2}(2\pi)^{m/2}} \frac{m!}{(m/2)!} \hat{V} ,
\end{equation}
where $P(W)$ is a 
polynomial in the Weyl curvature, and $\text{Pff}(W)$ is the Pfaffian evaluated in the Weyl curvature instead of the full curvature.
Equation \eqref{6d} shows that $P(W)$ is not zero in general.

As we mentioned in the introduction, Chang, Qing, and Yang \cite{Chang-Qing-Yang} have recently established a very similar formula \eqref{IntroCQY}.  It would be interesting to consolidate  these formulas.

\section{Varying the Poincar\'e-Einstein metric}\label{sec:Var}

In this section we will consider the variation of the characteristic forms of a Poincar\'e-Einstein manifold when we allow the metric to vary along a family $g_s$ of PE metrics with the same scalar curvature. 
After describing the expansion of $\dot{g}=\partial_s g_s |_{s=0}$, we recover results of Anderson \cite{Anderson} and Graham-Hirachi \cite{Graham-Hirachi} on the variation of the renormalized volume. Finally, we verify directly that the variation of the renormalized integral of the Pfaffian vanishes.

\subsection{Variation of the Pontrjagin Forms} 

As before, Pontrjagin forms and numbers on Poincar\'e-Einstein manifolds are easily dealt with due to conformal invariance.
Indeed, we know (Theorem \ref{Avez}) that
             \[ \tr\lrpar{ \bf{\Omega}^{\ell} }= \tr\lrpar{ \bar{\bf{\Omega}}^{\ell}} ,\]
so we need only consider variations of $\bar{g}$. Given a family of metrics $\bar{g}_s$ on $\bar{M}$
with connection and curvature forms $\bar{\omega}_s$, $\bar{\bf{\Omega}}_s$ and a homogeneous invariant polynomial in the curvature, $P(\bar{\mathbf{\Omega}}_s)$, of degree $q$ with complete polarization $p(\bar{\bf{\Omega}}_s,\ldots,\bar{\bf{\Omega}}_s)$, and
denoting the derivative at $s=0$ by an overdot, we have
\[ \dot{P}(\bar{\bf{\Omega}}) 
        = q \text{ d}p(\dot{\bar{\bf{\omega}}},\bar{\mathbf{\Omega}}, \ldots, \bar{\mathbf{\Omega}}) .\]

The derivative of a product of Pontrjagin forms is exact, being the product of closed and exact forms. Note that the boundary of $\bar{M}$ is totally geodesic, so restricting the Levi-Civita connections of $\bar{g}_s$ yields the Levi-Civita connections of the induced metrics. Hence $\dot{\bar{\omega}}$ when restricted to the boundary is the derivative of the connection 1-forms of the boundary metrics. We conclude with the following well-known result.

\begin{proposition}
Let $\bar{g}_s$ a smooth family of metrics on $\bar{M}$ such that, for every $s$, $(\bar{M}, \bar{g}_s)$ is a compact manifold with totally geodesic boundary. For any polynomial $Q$, the variation of the  characteristic number of $\bar{M}$ corresponding to $Q(\bf{\Omega})$ is the integral over the boundary of the Chern-Simons number corresponding to $Q(\bf{\Omega})$ on $\partial M$.
\end{proposition}

\subsection{Variation of the Renormalized Volume and the Pfaffian}\label{secVarPff}

Naturally, the renormalized Gauss-Bonnet theorem discussed previously shows that the variation of the renormalized integral of the Pfaffian vanishes. We verify this directly by studying the variation of the renormalized volume and the renormalized Weyl volume of tubes invariants.

Let $g_s$ be a family of Poincar\'e-Einstein metrics on $M$ with the same scalar curvature, $-m(m-1)$. We use $x$, a special boundary defining function for $g:=g_0$, to define a corresponding family of metrics on the boundary $g^{(0)}_s := x^2g_s \downharpoonright_{\partial M}$. The metrics $g^{(0)}_s$ in turn determine bdfs $x_s=xe^{\omega(x)}$, special with respect to $g_s$ (as in Lemma 2.1 of \cite{Graham}). We shall compute in a frame $X_u = x\bar{X}_u$ as in previous sections, and denote derivatives with respect to $s$ at $s=0$ by an overdot. Set $h=\dot{g}$.

Just as one needs to break conformal invariance by choosing a bdf in order to study the renormalized volume, the study of the Einstein equation requires breaking gauge invariance. The approach followed in \cite{Anderson} (see also \cite{Graham-Hirachi}) is to use the Fefferman-Graham expansion of the metrics $g_s$. A second approach is to impose the so called Bianchi gauge. One defines the operator
\[ B_g(k) := \delta_gk + \frac{1}{2}d\lrpar{ \tr_gh} \]
and demands that $g_s$ satisfy the Einstein equation and $B_g(g_s)=0$. This approach was first espoused in \cite{Biquard}, though we shall follow \cite{Mazzeo-Pacard} (see also \cite{Anderson2} and \cite{Lee}).

The following result appears as Lemma 3.4 in \cite{Anderson2}.

\begin{lemma}\label{hmm}
   $ x\partial_x \lrpar{\dot{\omega}} = \frac{1}{2}h_{mm}$.
\end{lemma}

\begin{proof}
We differentiate the relation $|dx_s|^2_{\bar{g}_s}=1$:
\begin{equation*}\begin{split}
            x_s^{-2} g_s \lrpar{ dx_s, dx_s }%
        &= x^{-2} 
             \lrspar{ g_s\lrpar{dx, dx} + 2xg_s\lrpar{ d\omega_s, dx} + %
             x^2 g_s \lrpar{d\omega_s, d\omega_s } } =1\\
\end{split}\end{equation*}
to get
\[ x^{-2} \lrspar{\dzero{s} g_s^{mm} + 2x g\lrpar{ dx, d\dot{\omega} } }  = 0 .\]
As $\partial g_s^{-1} = - g^{-1}\lrpar{ \partial g_s} g^{-1}$, we conclude that
\begin{equation*}
    x\partial_x \lrpar{\dot{\omega}} = \frac{1}{2}h_{mm}. 
\end{equation*}
\end{proof}

The linearized Einstein equation, $Ric'(h)=-(m-1)h$, is given by (1.179 in \cite{Besse})
\begin{equation}\label{LinEin}
            E_g'(h):= g^{st}(C_gh)_{uvs,t} -\frac{1}{2}Dd\tr_gh +(m-1)h =0 ,
\end{equation}
wherein
\[ (C_gh)_{uvs} = \frac{1}{2}\lrpar{ h_{sv,u} + h_{us,v} - h_{uv,s} } .\]
As the scalar curvature is the same for all $s$, we also have $Scal_g'(h)=0$, which by Theorem 1.174(e) of \cite{Besse}, means
\begin{equation}\label{scal} \begin{split}
     Scal_g'(h) 
     &=\Delta_g \lrpar{tr_g h} + \delta_g\lrpar{\delta_g h} - g\lrpar{Ric_g, g} \\
     &=\Delta_g \lrpar{tr_g h} + \delta_g\lrpar{\delta_g h} + (m-1) tr_g h = 0.
\end{split} \end{equation}

It is easy to compute the $0$-th order term in \eqref{LinEin}, for instance:
\begin{equation*}
         \lrpar{E_g'(h)}^{(0)}_{mm} = (m-1)h_{mm}^{(0)},
\end{equation*}
so $E_g'(h)=0$ implies $h_{mm}^{(0)}=0$, and we gain no more information from the other $0$-th order terms. Computing higher order terms will allow us to show that any perturbation of the metric through PE metrics has an expansion like that of the metric.

\begin{proposition}\label{hStructure}
Assume that $h_{mm}^{(0)} =0$, and that for some natural number $\alpha$ both:
\begin{equation}\label{Conditions} \begin{split}
 1)& \lrpar{X_u(h_{st})}^{(\alpha)} =  \delta_{um} \lrpar{ \alpha h_{st}^{(\alpha)} } \\
 2)& \lrpar{ \Gamma^{u}_{vw} h_{st}}^{(\alpha)} = %
                \lrpar{ \Gamma^{u}_{vw} }^{(\alpha)}\lrpar{ h_{st}}^{(0)}%
               +\lrpar{ \Gamma^{u}_{vw}}^{(0)}\lrpar{ h_{st}}^{(\alpha)}
\end{split}\end{equation}
hold, then
\begin{equation}\label{Ealpha}
     \lrpar{E_g'(h)}^{(\alpha)}_{uv} = \begin{cases}
            \frac{\alpha}{2} \lrpar{m-1-\alpha} h_{uv}^{(\alpha)} + \delta_{uv}%
               \lrspar{ \lrpar{m-1-\alpha} h_{mm}^{(\alpha)} + \frac{\alpha}{2} (\tr_gh)^{(\alpha)}}%
            & \text{ if $m \notin \{u,v\}$} \\
            \lrpar{1 - \frac{\alpha}{2}} %
            \lrspar{ \lrpar{m-1-\alpha} h_{mm}^{(\alpha)} + \alpha (\tr_gh)^{(\alpha)} }%
            & \text{ if $m=u=v$} \end{cases}.
\end{equation}

Consequently, once gauge invariance is broken, either through the Fefferman-Graham expansion or by imposing the Bianchi gauge, we have
\begin{equation}\label{hFeven}
           \mathcal{F} \lrpar{h_{uv}} = (-1)^{\delta_{um}+\delta_{vm}} ,
\end{equation}
($\mathcal{F}$ as in Corollary \ref{cor:Even} above), though in the odd-dimensional case log terms can occur with $x^{m-1}$.
\end{proposition}

\begin{proof}
The proof of \eqref{Ealpha} is a straightforward computation using \eqref{Conditions} in \eqref{LinEin}.

We first prove \eqref{hFeven} in the Fefferman-Graham expansion approach. We know that
\[ \bar{g}_s\lrpar{ \bar{\theta}^ m, \bar{\theta}^i} = 0 \text{ for any $i \neq m$.} \]
A computation like that of Lemma \ref{hmm} gives $h_{mi}= xg\lrpar{ d\dot{\omega}, \bar{\theta}^i}$.
This implies that $h_{mi}^{(0)}=0$ and, together with Lemma \ref{hmm}, that for $k<m-1$
\begin{equation}\label{hmi}
       \text{if $h_{mm}$ is even up to $x^k$, then $h_{mi}$ is odd up to $x^{k+1}$.}
\end{equation}
 This sets up an iterative scheme. In the first step, we use $h_{mi}^{(0)}=0$ to see that  the conditions \eqref{Conditions} are satisfied with $\alpha =1$. Then \eqref{Ealpha} shows that $h_{uv}^{(1)} = h_{mm}^{(1)}=0$, and this in turn, by \eqref{hmi}, that $h_{mi}^{(2)}=0$. Plugging this back into \eqref{Conditions} with $\alpha=3$ gives us the next iterative step. This works until $\alpha=m-1$, which establishes \eqref{hFeven}.

For the Bianchi gauge approach, consider the operator
\[ \Phi_g(k) = Ric_k + (m-1)k + \delta_k^*B_g(k) .\]
Clearly, $g_s$ is in the kernel of $\Phi_g$ and $h$ in that of $\Phi_g'$. We recall from Proposition 3 of \cite{Mazzeo-Pacard} that $\Phi_g'$ has no indicial roots between $0$ and $m-1$. This will yield \eqref{hFeven} through the following two observations.

First, we split $h$ into two pieces, $h=h_e + h_o$, by requiring that $(h_e)_{uv}$ have an even expansion in $x$ when there are an even number of $m$'s in $\{ u,v \}$ and have an odd expansion in $x$ otherwise.
It is not hard to see that $\Phi_g'(h_e)_{uv} \mod x^{m-1}$ is even or odd in $x$ when $(h_e)_{uv}$ is even or odd, and similarly $\Phi_g'(h_o)_{uv}$.

Secondly, $B_g(g_s)=0$ and $Scal'(h)=0$ (see \eqref{scal} above) together imply
that $\delta_g(h) =0$ and $tr_gh =0$.
An easy computation shows that $(\delta_gh)^{(0)}=0$ implies $h_{im}^{(0)}=0$.

These observations together show that $h_o$ vanishes to first order at the boundary, and that $\Phi_g'(h_o)$ vanishes to order $m-1$. We conclude that $h_o$ vanishes to order $m-1$ at the boundary. Thus $h = h_e \mod x^{m-1}$ and this is \eqref{hFeven}.
\end{proof}

We observe that \eqref{hFeven} shows that the conditions \eqref{Conditions} are satisfied with $\alpha=m-1$; this implies by \eqref{Ealpha} that $(tr_gh)^{(m-1)}=0$.

Once we know the structure of $h$, it is a simple matter to compute the variation of the renormalized volume in the even-dimensional case (Theorem 2.2 in \cite{Anderson}) and the variation of the residue in the odd-dimensional case (Theorem 1.1 in \cite{Graham-Hirachi}).

\begin{theorem}\label{V'}
On an even-dimensional PE manifold, if $g_s = g + sh$ to first order with $h$ as in Proposition \ref{hStructure}, the variation of the renormalized volume is
\begin{equation}
            \hat{V}'(h) = -\frac{1}{4} \int_{\partial M} \left< g^{(m-1)}, h^{(0)} \right> ,
\end{equation}
the inner product and the integral being taken with respect to $g^{(0)}$.
Similarly, on an odd-dimensional PE manifold, the variation of the residue of the volume is given by
\begin{equation}\label{L'}
         L'(h) = -\frac{1}{4} \int_{\partial M} \left< g^{(m-1,1)},h^{(0)} \right> .
\end{equation}
\end{theorem}

{\em Remark} In \cite{Graham-Hirachi}, \eqref{L'} is expressed in terms of Branson's $Q$ curvature and the Fefferman-Graham obstruction tensor $\curly{O}$ as
\[ \lrpar{ \int_{\partial M} Q }'(h) 
      = (-1)^{\frac{m-1}{2}} \frac{m-3}{2} \int_{\partial M} \left< \curly{O}, h \right> .\]

\begin{proof}
Using Riesz renormalization, we have
\begin{equation}\label{Eq1}   \begin{split}
           \dzero{s} \FP_{z=0} \int_{M} x_s^z \dvol_s  %
           &= \FP_{z=0} \int_M zx^{z-1} \dot{x} \dvol %
                 + \FP_{z=0} \int_M x^z \lrpar{ - \frac{1}{2} tr_g h} \dvol \\
          &= \int_{\partial M} \dot{\omega}^{(m-1)} \dvol_{0}  %
               + \FP_{z=0} \frac{1}{2(m-1)} \int_M x^z \lrspar{ \Delta \lrpar{tr_g h}%
               + \delta_g\delta_g h} \dvol \\
          &= \int_{\partial M} \lrspar{ %
                \dot{\omega}^{(m-1)} + \frac{1}{2} \lrpar{ tr_g h}^{(m-1)} 
               + \frac{1}{2(m-1)}\lrpar{\delta_g h}_{m}^{(m-1)}  } \dvol_{0}\\  %
          &=  \frac{1}{2(m-1)}\int_{\partial M} \lrspar{ %
                h_{mm}^{(m-1)} +\lrpar{\delta_g h}_{m}^{(m-1)}  } \dvol_{0}.  %
\end{split} \end{equation}
The second equality is \eqref{scal}, the third follows from integrating by parts,
and the last uses Lemma \ref{hmm} and the observation made before the statement of the theorem.

Now as in Proposition \ref{hStructure}, we can use the structure of $h$ to compute
\begin{equation}\label{Eq2}
       \lrpar{\delta_g h}_{m}^{(m-1)} = -h_{mm}^{(m-1)} - \frac{m-1}{2}  \left< g^{(m-1)}, h^{(0)} \right>  ,
\end{equation}
which gives the theorem in the even-dimensional case.

In odd dimensions the same proof works. The reason is that the $x^{m-1}\log x$, being the first $log$ term, behaves like the first odd term in the even-dimensional case. Specifically, it satisfies conditions \eqref{Conditions} with $\alpha=(m-1,1)$ which allows us to compute as in \eqref{Eq1} and \eqref{Eq2} replacing $(m-1)$ by $(m-1,1)$.
\end{proof}

To deal with the other terms in the expression for the Pfaffian in section $\mathcal{x}$\ref{sec:Pff} we 
use the formalism of double forms, e.g., we interpret $h$ as an element of $\curly{D}^{1,1}$ and denote it by $\df{h}$. 
We benefit greatly from a recent description of the variation of Weyl's volume of tubes invariants by Labbi \cite{Labbi}, \cite{Labbi0}. We denote the integrand in the definition of $K_{2\ell}$ by
\[ k_{2\ell}(\df{R}) := \frac{\curly{C}^{2\ell}\df{R}^{\ell}}{(2\ell)!(\ell)!} .\]
Labbi defines the tensors
\begin{equation}\label{GenEin}
       \df{E}_{2\ell}(\df{R}) := k_{2\ell}(\df{R}) \cdot \df{g} 
                                         - \frac{ \curly{C}^{2\ell -1}\df{R}^{\ell}}{(2\ell -1)! \ell !} ,
\end{equation}
and establishes that, on any closed Riemannian manifold,
\begin{equation}\label{KR'}
      \lrpar{K_{2\ell}(\df{R})}'(\df{h}) = \int \left< \df{E}_{2\ell}(\df{R}), \df{h} \right> .
\end{equation}
We will show that a very similar formula holds, on Einstein manifolds, for the variation of the Pfaffian. Our definitions consistently differ from those of Labbi by a factor of $\ell!$ but agree with those of \cite{Tubes}.

The tensors \eqref{GenEin} were first introduced by Lovelock in \cite{Lovelock}. His interest was to find all divergence-free symmetric (0,2)-tensors built from the metric and its first two covariant derivatives. His theorem is that they are given by arbitrary linear combinations of the $\df{E}_{2\ell}(\df{R})$ and the metric. These properties are  sometimes used to motivate the Einstein field equation in four dimensions (e.g., 3.7(i) in \cite{Besse}). Lovelock's generalization plays an analogous role in higher dimensional GB gravity theories (see \cite{Padilla} and references therein). In the physics literature, a linear combination of the $k_{2\ell}(\df{R})$ is known as a Lovelock Lagrangian; an arbitrary linear combination of the $\df{E}_{2\ell}(\df{R})$ is known as a Lovelock tensor. Note that $\df{E}_2(\df{R})$ is the usual stress energy tensor, $\frac{Scal}{2}g-Ric$. 
 Recall that any surface is automatically Einstein and that any higher dimensional Einstein manifold has constant scalar curvature. The following results appear in \cite{Patterson}.

\begin{lemma}\label{EmW=0}
      On any $n$ dimensional manifold, \newline
     a) {\em (Bach-Lanczos identity)} If $n$ is even,
            \[ \df{E}_n(\df{R}) = \df{E}_n(\df{W}) = 0.\]
     b) If $2\ell <n$ and $\df{E}_{2\ell}(\df{R})= \lambda\df{g}$ for some constant $\lambda$, then $k_{2\ell}(\df{R})$ is constant.
\end{lemma}

\begin{proof} $ $\newline
a) Any double form in $\curly{D}^{n,n}$ is of the form $f \df{g}^{n}$ for some function $f$, and it is easy to see that
\[ \frac{\curly{C}^{n}(f \df{g}^n)}{n! \lrpar{\frac{n}{2}}!} \df{g} 
             - \frac{\curly{C}^{n-1}(f \df{g}^{n})}{(n-1)! \lrpar{\frac{n}{2}}!} = 0. \]
b) We always have
\[ \curly{C}\lrpar{\df{E}_{2\ell}(\df{R})} 
         = \frac{\lrpar{n-2\ell}}{ \lrpar{2\ell}! \ell !} \curly{C}^{2\ell}\df{R}^{\ell}, \]
so if in addition we know that $\curly{C}\lrpar{\df{E}_{2\ell}(\df{R})}= n\lambda$ and $n \neq 2\ell$, we find that $\curly{C}^{2\ell}\df{R}^{\ell}$ is constant.
\end{proof}

Labbi's proof of \eqref{KR'} follows from this formula, Lemma 4.2 of \cite{Labbi} ($\ell \geq 1$):
\begin{equation}
        \lrpar{k_{2\ell}(\df{R})}'(\df{h}) = 
               -\frac{1}{2}  \left< \frac{ \curly{C}^{2\ell-1}\df{R}^{\ell} }{(2\ell -1)!(\ell !)}, \df{h} \right>
               +  \frac{(-1)^m}{4} \lrpar{\delta \widetilde{\delta} + \widetilde{\delta}\delta}
                   \lrpar{ * \lrpar{ \frac{ \df{g}^{m-2\ell} \cdot \df{R}^{\ell-1} }{ (m-2\ell)! (\ell -1)!}  \cdot \df{h}}} .
\end{equation}
where $*$ extends to double forms via
\[ * \lrpar{ \omega \otimes \eta } = *\omega \otimes *\eta ,\]
and 
 $\lrpar{\delta \widetilde{\delta} + \widetilde{\delta}\delta}: \curly{D}^{1,1} \to \curly{D}^{0,0}$ is
 a second order differential operator; all we will need to know here is that its adjoint is twice the Hessian (see comments after equation (10) in \cite{Labbi}), so that
 \[ \lrpar{\delta \widetilde{\delta} + \widetilde{\delta}\delta}^{*}(x^z) = 2z(-1)^m  x^z \FunForm ,\]
with $\FunForm \in \curly{D}^{1,1}$ the second fundamental form of $g$.

For a variation of Einstein metrics preserving the scalar curvature, a very similar formula holds.
The following proposition is proved by adapting Labbi's approach from $\df{R}$ to $\df{W}$.

\begin{proposition}\label{thm:varPff}
Consider a family $g_s$ of PE metrics on an even-dimensional manifold, with $g_s = g + sh$ to first order with $h$ as in Proposition \ref{hStructure}. \newline
a) For $\ell>2$
 \[ \dzero{s} \lrpar{ \Rint k_{2\ell}(\df{W})} = 
                \Rint \lrspar{
                 \frac{\left< \df{E}_{2\ell}(\df{W}) , \df{h} \right>}{2}  %
                  - \lrpar{\frac{s_g}{2m(m-1)}} \lrpar{m-2\ell +1} \left< \df{E}_{2\ell-2}(\df{W}), \df{h} \right> }, \]

and for $\ell =2$,
\[ \dzero{s} \lrpar{ \Rint k_{4}(\df{W})} =
                \frac{1}{2} \Rint \left< \df{E}_{4}(\df{W}) , \df{h} \right> 
                +\frac{(m-1)(m-3)}{4} \int_{\partial M}  \left< g^{(m-1)}, h^{(0)} \right>_{g^{(0)}}.\]
b) For these variations, the variation of the renormalized integral of the Pfaffian vanishes.
\end{proposition}

{\em Remark.}
Notice that $\df{E}_2(\df{W}) =0$ because $\curly{C}\df{W}=0$. In four dimensions, the functional
\[ g \mapsto K_4(\df{W}) = \int |W|^2 \]
is well understood. Its gradient is the Bach tensor which vanishes for metrics conformal to Einstein metrics. This is reflected in the vanishing of the interior integral in the formula above for $\ell=2$, since by the Bach-Lanczos formula $\df{E}_4(\df{W})=0$ in four dimensions.
It would be interesting to understand the behavior of the $K_{2\ell}(\df{W})$ under an arbitrary variation of the metric.

\begin{proof}
Define
\begin{equation*}\begin{split} 
       F_{\df{h}}(\df{R})&\lrpar{(X,Y)(Z,W)} \\
       &= h\lrpar{R(X,Y)Z,W} -h\lrpar{R(X,Y)W,Z} +h\lrpar{R(Z,W)X,Y} - h\lrpar{R(Z,W),Y,X} .
\end{split} \end{equation*}
Lemma 4.1 of \cite{Labbi} establishes
\[ \df{R}'(\df{h}) 
           = -\frac{1}{4}\lrpar{D\widetilde{D} + \widetilde{D}D}(\df{h}) + \frac{1}{4}F_{\df{h}}(\df{R}), \]
for some differential operators $D, \widetilde{D}$.
Simple manipulations using that the metrics $g_s$ all have the same scalar curvature yield
\begin{equation}\label{W'h}\begin{split}
      \df{W}'(h) &= \df{R}'(\df{h}) - \frac{s}{2m(m-1)} \lrpar{\df{g}^2}'(\df{h})\\
      &= -\frac{1}{4}\lrpar{D\widetilde{D} + \widetilde{D}D}(\df{h}) + \frac{1}{4}F_{\df{h}}(\df{R}) - 
              \frac{2s}{2m(m-1)} \lrpar{\df{g}\cdot \df{h}}\\
      &= -\frac{1}{4}\lrpar{D\widetilde{D} + \widetilde{D}D}(\df{h}) + \frac{1}{4}F_{\df{h}}(\df{W}) - 
              \frac{s}{2m(m-1)} \lrpar{\df{g}\cdot \df{h}}
\end{split}\end{equation}
Labbi's computations in the proof of Lemma 4.2 of \cite{Labbi} go through with $\df{R}$ replaced with $\df{W}$ and $\df{R}'(\df{h})$ replaced by \eqref{W'h} because, on Einstein manifolds, Weyl curvature satisfies the second Bianchi identity. This gives
\begin{equation}\label{vark} \begin{split}
        \lrpar{k_{2\ell}(\df{W})}'(\df{h}) &= %
               -\frac{1}{2}  \left< \frac{ \curly{C}^{2\ell-1}\df{W}^{\ell} }{(2\ell -1)!(\ell !)}, \df{h} \right>%
               - \lrpar{\frac{s_g}{2m(m-1)}} \lrpar{m-2\ell +1} \left< \df{E}_{2(\ell-1)}(\df{W}), \df{h} \right> \\
               &+  \frac{(-1)^m}{4} \lrpar{\delta \widetilde{\delta} + \widetilde{\delta}\delta}%
                   \lrpar{ * \lrpar{ \frac{ \df{g}^{m-2\ell} \cdot \df{W}^{\ell-1} }{ (m-2\ell)! (\ell -1)!}  \cdot \df{h}}} .
\end{split} \end{equation}

Consider the renormalized integral of the last term in equation \eqref{vark},
\begin{equation*}\begin{split}
  \FP_{z=0} \int_M x^{z} (-1)^m \lrpar{\delta \widetilde{\delta} + \widetilde{\delta}\delta}
                   \lrpar{ * \lrpar{ \df{g}^{m-2\ell} \cdot \df{W}^{\ell-1} \cdot \df{h}}} 
     &= 2\FP_{z=0} \int_M z  x^z \left<  \FunForm, 
                   * \lrpar{ \df{g}^{m-2\ell} \cdot \df{W}^{\ell-1} \cdot \df{h}} \right> \\
     &= 2 \Res \lrspar{ \left< \FunForm, * \lrpar{ \df{g}^{m-2\ell} \cdot \df{W}^{\ell-1} \cdot \df{h}} \right>}. \\
\end{split}\end{equation*}
So this term produces residues, and we need to find the  $x^{m-1}$ term in the expansions. We revisit Lemma \ref{lem:Gammas} and ascertain that
\begin{equation*}
         \left< R(X_i,X_j)X_k, X_\ell \right>_g ^{(0)} =
              \lrspar{ \gamma^m_{ik}\gamma^m_{j\ell} - \gamma^m_{jk}\gamma^m_{i\ell} }^{(0)} ,
\end{equation*}
and 
\begin{equation*}
         \left< R(X_i,X_j)X_k, X_\ell \right>_g ^{(m-1)} =
               \lrspar{ \gamma^m_{ik}\gamma^m_{j\ell} - \gamma^m_{jk}\gamma^m_{i\ell} }^{(m-1)} .
\end{equation*}
Thus the coefficients of $\df{W}$ computed in this frame have no constant term, and the residue vanishes if $\ell >2$. To compute the residue  for $\ell =2$, denote the double form corresponding to $g^{(0)}$ by $\df{g}_{0}$ and similarly $\df{g}_{m-1}$. Note that replacing the coefficients of $\FunForm$ by their $0$-th order part produces the double form $\df{g}_0$ and replacing the coefficients of $\df{W}$ by their $(m-1)$-th order part
produces the double form $\df{g}_{0} \cdot \lrpar{ -\frac{m-1}{2}\df{g}_{m-1} }.$
We can use the following formulas from \cite{Labbi0}
\begin{equation*}
      \left< \omega, * \eta \right> = \left< * \omega, \eta \right> \text{, }
      \left< \df{g} \cdot \omega, \eta \right> = \left< \omega, \curly{C} \eta \right> \text{, }
      *\frac{\df{g}^k}{k!} = \frac{\df{g}^{m-k}}{(m-k)!} \text{, and }
       \curly{C}^k\lrpar{\frac{\df{g}^\ell}{\ell !} }= \frac{(m-\ell +k)!}{(m-\ell)!}\frac{\df{g}^{(\ell -k)}}{(\ell - k)!},
\end{equation*}
to see that
\begin{equation}\label{k=4} \begin{split}
         \frac{1}{2} \Res & \left[ \left< \FunForm  , %
                          * \lrpar{ \frac{\df{g}^{m-4} \cdot \df{W}}{(m-4)!} \cdot \df{h}} \right> \right] %
        =-\frac{m-1}{4(m-4)!} \int_{\partial M} %
                 \left< \df{g}_0, * \lrpar{\df{g}_0^{m-4} \cdot \df{g}_0 \cdot \df{g}_{m-1} \cdot \df{h}_0 } \right>\\
        &\phantom{xxxxxxxxxxx}=-\frac{m-1}{4(m-4)!} \int_{\partial M} %
                  \left< (m-3)! \frac{\df{g}^2}{2} , \df{g}_{m-1} \cdot \df{h}_0 \right>\\
        &\phantom{xxxxxxxxxxx}=-\frac{(m-1)(m-3)}{4} \int_{\partial M} %
                  \lrspar{ \lrpar{\tr_{g^{(0)}} g^{(m-1)}} \lrpar{\tr_{g^{(0)}} h^{(0)}} %
                          - \left< g^{(m-1)}, h^{(0)} \right>_{g^{(0)}} }\\
        &\phantom{xxxxxxxxxxx}=\frac{(m-1)(m-3)}{4} \int_{\partial M} %
                  \left< g^{(m-1)}, h^{(0)} \right>_{g^{(0)}}.
\end{split} \end{equation}

Putting these observations together, we see that the variation of the $\ell$-th renormalized volume of tube invariant evaluated in the Weyl curvature
       \[ K_{2\ell}(\df{W}) =  \Rint \frac{\curly{C}^{2\ell}(\df{W}^{\ell})}{(2\ell)!\ell ! } \dvol \]
is given by
\begin{equation*}   \begin{split}
           \lrpar{K_{2\ell}(\df{W})}'(\df{h})  %
           &= \FP_{z=0} \int_M zx^{z} \dot{\omega} k_{2\ell} \dvol %
                 + \FP_{z=0} \int_M x^z \dot{k}_{2\ell} \dvol %
                 + \FP_{z=0} \int_M x^z k_{2\ell} \lrpar{ \frac{1}{2} tr_g h} \dvol \\
           &=  \FP_{z=0} \int_M x^z \lrspar{%
                 \frac{\left< \df{E}_{2\ell}(\df{W}) , \df{h} \right>}{2}  %
                  - \lrpar{\frac{s_g}{2m(m-1)}} \lrpar{m-2\ell +1} \left< \df{E}_{2\ell-2}(\df{W}), \df{h} \right> }
                  \dvol \\
           &\phantom{xx}      +\Res \lrpar{ \dot{\omega} k_{2\ell} }
                 +\frac{1}{2} \Res \lrspar{ \left< \FunForm, %
                   * \lrpar{ \frac{ \df{g}^{m-2\ell} \cdot \df{W}^{\ell-1}}{(m-2\ell)!(\ell-1)!} \cdot \df{h}} \right>}. \\
\end{split} \end{equation*}
The second residue has been computed above for $\ell=2$ and vanishes otherwise.
For the first residue, we know that $\df{W} = x^2 \bar{\df{W}}$ (e.g., \cite{Kulkarni}) hence neither $\dot{\omega}$ nor $k_{2\ell}(\df{W})$ have a constant term, and their product does not have an $x^{m-1}$ term in its expansion, so that residue also vanishes. This establishes part (a) of the proposition.

Using  the expression in Lemma \ref{EinPfaff} for the Pfaffian of an Einstein metric, its variation can be read off from part (a) of the proposition and the variation of the volume in Theorem \ref{V'}. The interior integrals telescope
\begin{equation*}\begin{split}
          \frac{1}{(2\pi)^{m/2}}&  \Rint \sum_{k=0}^{m/2-2} \frac{(2k)!}{k!} \tilde{s}^k
             \lrpar{ \frac{\left< \df{E}_{m-2k}(\df{W}) , \df{h} \right>}{2}  %
                  - \tilde{s} \lrpar{2k +1} \left< \df{E}_{m-2k-2}(\df{W}), \df{h} \right> } \\
          =\frac{1}{(2\pi)^{m/2}}&  \Rint \sum_{k=0}^{m/2-2} \frac{(2k)!}{k!} \tilde{s}^k
             \lrpar{ \frac{\left< \df{E}_{m-2k}(\df{W}) , \df{h} \right>}{2}  %
                  - \tilde{s} \frac{(2k +1)(2k+2)}{(k+1)}%
                  \frac{\left< \df{E}_{m-2k-2}(\df{W}), \df{h} \right>}{2} }  ,
 \end{split}\end{equation*}
where we have abbreviated the scalar curvature factor to $\tilde{s}$. We are left only the $\left< \df{E}_m(\df{W}),\df{h} \right>$ term, which is zero by Lemma \ref{EmW=0}. There are two terms with residues, one coming from the volume and the other from $\ell =2$. The corresponding summands in the formula for the Pfaffian (in Lemma \ref{EinPfaff}) are
\[ \frac{ m!}{(\frac{m}{2})!} \frac{(-1)^{\frac{m}{2}}}{2^{\frac{m}{2}}}
                  \lrspar{ 1 + \frac{k_4}{(m-1)(m-3)} } ,\]
so from \eqref{k=4} and Theorem \ref{V'}, we see that these residues cancel each other out.
\end{proof}


\end{document}